\documentclass[11pt]{amsart}
 \usepackage{amsopn}
 \usepackage{amsmath,amsthm,amssymb}

 \newcommand{\nc}{\newcommand}
 
 \nc{\ag}{\mathfrak{a} } \nc{\bb}{\mathfrak{b} }
 \nc{\cc}{\mathfrak{c} }  \nc{\dd}{\mathfrak{d} } 
 \nc{\ggo}{\mathfrak{g} }
 \nc{\hh}{\mathfrak{h} }  \nc{\ii}{\mathfrak{i} }
 \nc{\jj}{\mathfrak{j} }  \nc{\kk}{\mathfrak{k} }
\nc{\mm}{\mathfrak{m} }   \nc{\nn}{\mathfrak{n} }
\nc{\pp}{\mathfrak{p} }   \nc{\sg}{\mathfrak{s} }
 \nc{\sog}{\mathfrak{so} }  \nc{\spg}{\mathfrak{sp} }
 \nc{\sug}{\mathfrak{su} }  \nc{\slg}{\mathfrak{sl} }
 \nc{\tg}{\mathfrak{t} }  \nc{\uu}{\mathfrak{u} }
 \nc{\vv}{\mathfrak{v} } \nc{\ww}{\mathfrak{w} }
 \nc{\zz}{\mathfrak{z} }  
 
 \nc{\ggob}{\overline{\mathfrak{g}}} 
 
\nc{\glg}{\mathfrak{gl} }
  
\nc{\pca}{\mathcal{P}} \nc{\nca}{\mathcal{N}}
 
 \nc{\vp}{\varphi} \nc{\ddt}{\frac{{\rm d}}{{\rm d}t}}
 \nc{\la}{\langle} \nc{\ra}{\rangle}
 
 \nc{\SO}{{\sf SO}} \nc{\Spe}{{\sf Sp}} \nc{\Sl}{{\sf Sl}}
 \nc{\SU}{{\sf SU}} \nc{\Or}{{\sf O}} \nc{\U}{{\sf U}}
 \nc{\Gl}{{\sf Gl}} \nc{\Se}{{\sf S}} \nc{\Cl}{{\sf Cl}}
 \nc{\Spin}{{\sf Spin}} \nc{\Pin}{{\sf Pin}}
 
 \nc{\RR}{{\mathbb R}} \nc{\HH}{{\mathbb H}} \nc{\CC}{{\mathbb C}}
 \nc{\ZZ}{{\mathbb Z}} \nc{\FF}{{\mathbb F}} \nc{\NN}{{\mathbb N}}
 \nc{\GG}{{\mathbb G}} \nc{\JJ}{{\mathbb J}} \nc{\II}{{\mathbb I}}
 \nc{\KK}{{\mathbb K}} \nc{\DD}{{\mathbb D}}
 
 \nc{\ad}{\operatorname{ad}} \nc{\Ad}{\operatorname{Ad}}
 \nc{\coad}{\operatorname{coad}}
 \nc{\rank}{\operatorname{rank}} \nc{\Irr}{\operatorname{Irr}}
 \nc{\End}{\operatorname{End}} \nc{\Aut}{\operatorname{Aut}}
 \nc{\Inn}{\operatorname{Inn}} \nc{\Der}{\operatorname{Der}}
 \nc{\Ker}{\operatorname{Ker}} \nc{\Iso}{\operatorname{I}}
 \nc{\Le}{\operatorname{L}} \nc{\tr}{\operatorname{tr}}
 \nc{\dif}{\operatorname{d}} \nc{\sen}{\operatorname{sen}}
 \nc{\modu}{\operatorname{mod}} \nc{\Ric}{\operatorname{R}}
 \nc{\Sym}{\operatorname{Sym}} \nc{\sca}{\operatorname{sc}}
 \nc{\scalar}{{\sf s}} \nc{\grad}{\operatorname{grad}}
 \nc{\ricci}{\operatorname{r}} \nc{\riccin}{\operatorname{Ric}}
 \nc{\Lie}{\operatorname{L}} \nc{\tang}{\operatorname{T}}
 
 \theoremstyle{plain}
 \newtheorem{thm}{Theorem}[section]
 \newtheorem{prop}[thm]{Proposition}
 \newtheorem{cor}[thm]{Corollary}
 \newtheorem{lem}[thm]{Lemma}
 
 \theoremstyle{definition}

 \theoremstyle{remark}
 \newtheorem*{rem}{Remark}
 
 \newtheorem{exam}[thm]{Example}
 \newtheorem{exams}[thm]{Examples}

 \newcommand{\ri}{{\rm (i)}}
 \newcommand{\rii}{{\rm (ii)}}
 \newcommand{\riii}{{\rm (iii)}}

\begin{document}
 \title[Two-step nilpotent Lie algebras with ad-invariant metrics] {two-step nilpotent Lie algebras with ad-invariant metrics and a special kind of skew-symmetric maps}
\author[Gabriela Ovando]{Gabriela Ovando}
 
 \address{Gabriela Ovando, Dpto. de Matem\'atica, ECEN - FCEIA, Pellegrini 250, 2000 Rosario, Santa Fe, Argentina.}
 
 \address{FaMAF, Ciudad Universitaria, 5000 C\'ordoba Argentina}                                    
    %\curraddr{FCEIA}                                   
    \email{ovando@mate.uncor.edu}                                      
    %\urladdr{...www}                                    
    %\dedicatory{to my sun}                                 
    %\date{November 2006}                                       
     
    \thanks{Partially supported by CONICET, FONCYT and SECYT (UNC)}                                               
    %\translator{martin}                                 
    %\subjclass{22E25, 22E60, 53B30}                                  
    \begin{abstract} We prove that a  2-step nilpotent Lie algebras admitting an ad-invariant metric can be constructed from a vector space $\vv$ endowed with a inner product $\la \,,\,\ra$ and an injective homomorphism $\rho: \vv \to \sog(\vv)$ satisfying $\rho(v)v=0$ for all $v\in \vv$. The corresponding simply connected pseudo-Riemannian Lie groups are flat and  any isometry fixing the identity element does not depend on $\rho$. The description allows to construct examples starting with a compact semisimple Lie algebra and it is useful to show some  applications. \end{abstract}
 
 \maketitle
 
 \section{Introduction}

  Two-step nilpotent Lie groups endowed with a left invariant Riemannian metric have been considered with special attention in Riemannian geometry, harmonic analysis and spectral geometry. But this Riemannian metric cannot be also right invariant. If we ask a smoothly non degenerate symmetric form on a 2-step nilpotent Lie group to be left and right invariant then the form must have a non trivial index unless the group is abelian. The algebraic counterpart of a Lie group endowed with a bi-invariant metric is a 2-step nilpotent Lie algebra provided with an ad-invariant metric. As we shall see this condition of the metric to be ad-invariant determines the algebraic structure of the Lie algebra.  An open problem is whether a given Lie algebra can be provided with an ad-invariant metric and the attempts to get a complete answer to this question promote mathematicians to study the structure of a Lie algebra provided with an ad-invariant metric. For instance a  known necessary condition for a solvable Lie   algebra $\ggo$ to admit such metric is the following one: $\dim \ggo= \dim \zz(\ggo) + \dim C^1(\ggo)$, where $\zz(\ggo)$ and $C^1(\ggo)$ denote the center and commutator of $\ggo$ respectively. In case of being $\ggo$ 2-step nilpotent one gets $\dim C^1(\ggo) = \dim (\ggo / \zz(\ggo))$. However this is not sufficient. 
 Lie algebras with an ad-invariant metric were described as the result of a double extension procedure by different authors (see \cite{Me} \cite{MR} \cite{FS}). This can be done by a one to one step construction but such that from one step to the next one the algebraic structure of the resulting Lie algebra is quite different from the  original ones.  In the case of 2-step nilpotent Lie algebras admitting an ad-invaraint metric,  Noui and Revoy \cite{NR} follows a different approach by relating this kind of Lie algebras with alternating trilinear forms.
  
  Our aim in this work is to give a description of the 2-step nilpotent Lie algebras provided with an ad-invariant metric following a setting comparable to that used by people working with nilmanifolds provided with a left invariant Riemannian metric. Let us recall this.  
 In the  80's Kaplan produced a 2-step nilpotent Lie group with a left invariant metric whose Lie algebra arises  from a real representation $(J,V)$ of a Clifford algebra $Cl(\zz)$. This construction was generalized in the following years showing that any 2-step nilpotent Lie algebra endowed with an inner product can be constructed starting with a vector space $(\vv, \la \, , \, \ra_{\vv})$, and a linear homomorphism $J: \zz \to \sog(\vv)$. Furthermore,  it was proved that many geometric and analytical properties of the respective Lie groups  can be read in terms of $J$.    
 
   We prove that 2-step nilpotent  Lie algebras provided with an ad-invariant metric can be achieved from a injective homomorphism $\rho:\vv \to \sog(\vv)$ of a real vector space $\vv$ equipped with an inner product $\la \, , \, \ra$  satisfying
 $$ (*) \qquad \qquad \rho(v)v=0 \qquad\mbox{ for all }\quad v\in \vv. $$
  
   A  consequence of this result is that  no 2-step nilpotent Lie algebra with one, two or four dimensional commutator  can be endowed with an ad-invariant metric, while for $\dim \vv$ equals  3 or greather than four, there always exists a Lie algebra admitting an ad-invariant metric, whose commutator is isomorphic as vector space with $\vv$. Furthermore the existence of an   ad-invariant metric  on a   2-step nilpotent Lie algebra implies that $\nn$  cannot be non singular.  In particular H-types Lie algebras cannot be provided with an ad-invariant metric.
 
 With our description we construct examples and we show applications of them.  In this sense our study   allows to get a 2-step nilpotent Lie algebra $\nn$ which can be attached with an ad-invariant metric  starting  with a semisimple compact Lie algebra $\ggo$ and a representation of it.  We  applied  this to give double Lie algebras $(\ggo, r)$ where $r$ is a classical r-matrix. Up an $r$-matrix on a compact semisimple Lie algebra $\ggo$ we are able to construct an $r$-matrix on the corresponding 2-step nilpotent Lie algebra $\nn$ getting from $\ggo$. Furthermore if $(\ggo,r)$ gives rise to a Lie bialgebra structure then $\nn$ adquires a  Lie bialgebra structure. Finally we point out some features about the geometry derived from the bi-invariant metric. In the simply connected case these Lie groups are flat and the group of isometries fixing the identity element does not depend on the Lie bracket, hence  the isometry group can be much larger than the group of isometric authomorphisms, a property that makes a difference with the Riemannian case.
 
 \section{Preliminaries}
  In this section we give  basic definitions and we expose some features to be used through the paper. We assume all Lie algebras are real. 
  
 \subsection{Vector spaces with (non definite) metrics}  A metric  on a real vector space $\vv$ is a symmetric bilinear form on $\vv$, $\la \, , \,  \ra:\vv \times \vv \to \RR$ which is non degenerate, that is, for any non zero vector $x\in \vv$ there exists a vector $y\in\vv$ such that $\la x, y\ra \neq 0$. If this does not occur the symmetric bilinear form $\la \, , \, \ra$ is said degenerate.
   
   If $\ww$ is a subspace of  $(\vv, \la \, , \, \ra)$ the subspace
 $$\ww^{\perp} = \{ x \in \vv \,:\, \la x, v\ra =0 \quad \text{ for all } \quad v \in \ww\}$$
 denotes the {\it orthogonal}  subspace of $\ww$. In particular we say that $\ww$ is {\em isotropic} if 
 $\ww \subset \ww^{\perp}$, {\em totally isotropic} if $\ww =\ww^{\perp}$  and {\em non degenerate} if $\ww \cap \ww^{\perp}=0$. In the last case the restricted metric $\la \, , \, \ra_{\ww}:=\la \, , \, \ra_{|_{\ww \times \ww}}$ defines a isomorphism $\xi$ between $\ww$ and its dual space $\ww^{\ast}$ by $\xi(u)(v)=\la u,v\ra_{\ww}$. As usual a metric of index 0 is called a {\it inner product}.
 
 The proof of the following lemma follows by an induction procedure on the dimension of $\vv$ by application of the above definitions.
 
 \begin{lem}\label{l1} Let $(\vv, \la \,, \,\ra)$ be a vector space equipped with a metric and assume that $\uu \subset \vv$ is  a totally isotropic subspace of $\vv$. Then there exists a complement $\ww$ of $\uu$ in $\vv$ such that
 $$\vv= \uu \oplus \ww$$ 
 where $\ww$ can be choosen either as a totally isotropic subspace or either as a non degenerate subspace on which the retricted metric is a inner product on $\ww$.
 \end{lem}
 
 Notice that the dimension of the vector space $\vv$ in the previous lemma must be even and exactly $2 \dim \uu$. The index of the metric is therefore $\dim \uu = \frac12 \dim \vv$. Such metric  is called {\em hyperbolic} or {\em neutral}. Conversely let  $\uu$ denote a vector space  whose dual space is  denoted by $\uu^{\ast}$. Let $\uu^*\oplus \uu$ be the direct sum as vector spaces of $\uu$ and $\uu^*$ and endow this with the hyperbolic metric $\la \phi_1 + x_1, \phi_2 + x_2 \ra = \phi_1 (x_2) + \phi_2(x_1)$ where $\phi_i \in \uu^*$, $x_i\in \uu$, for  i=1,2. Clearly  $\uu^*$ and $\uu$ are complementary totally isotropic subspaces in $\uu^* \oplus \uu$.

If $(\vv, \la \, ,\, \ra)$ is a vector space endowed with a metric $\la \,, \,\ra$, then $\sog(\vv)$ denotes the set of skew-symmetric linear transformations in $\vv$.

 %\vskip 5pt
 
 \subsection{On Lie algebras provided with an ad-invariant metric} Let $\ggo$ denote a Lie algebra. A symmetric bilinear form  $\la \,, \,\ra$ on $\ggo$ is called {\em ad-invariant} if it satisfies
 \begin{equation}\label{ad}
 \la [x,y], z \ra + \la y,[x,z]\ra=0 \quad \text{ for all } \quad x,y, z \in \ggo.
 \end{equation}
 Thus the image of the coadjoint representation $\ad:\ggo \to \End(\ggo)$, $\ad(x) (y) = [x,y]$  belongs to $\sog(\ggo)$.
 
 Easy computations show that any ad-invariant symmetric bilinear form on the three dimensional Heinsenberg Lie algebra $\hh_3$ is degenerate, where $\hh_3=span\{x,y,z\}$ with the Lie bracket given by $[x,y]=z$ and the other relations are zero. This says that not every Lie algebra can be provided with an ad-invariant metric.
 
 \begin{exams} \label{cot} The Killing form $B$,  $B(x,y)=tr(\ad_x \ad_y)$  where $tr$ denotes the trace, provides an ad-invariante metric on any semisimple Lie algebra. 
 %In particular if $\ggo$ is also compact, then $-B$ defines an ad-invariant inner product on the semisimple Lie algebra  $\ggo$.
 
 Examples of Lie algebras with an ad-invariant metric can be obtained in the following way. Let $\hh$ denote a Lie algebra, whose dual vector space  $\hh^*$ is an abelian  Lie algebra with the trivial Lie bracket. This is an $\hh$-module with the {\em coadjoint action}
 $$x \cdot \phi = -\phi \circ \ad_x \quad \text{ for }x \in \hh,  \phi \in \hh^*.$$
  The semidirect product $\hh^* \rtimes \hh$, called the {\em cotangent of $\hh$}, can be equipped with an ad-invariant metric if we put
 $$\la \phi_1 + x_1, \phi_2 + x_2 \ra = \phi_1 (x_2) + \phi_2(x_1) \qquad \phi_i \in \hh^*, \, x_i\in \hh,\, i=1,2$$ 
 which is the canonical hyperbolic metric on $\hh^*\rtimes\hh$. 
 \end{exams}
 
 If we ask an  $\ad$-invariant metric on a real Lie algebra to be an inner product then $\ggo$ is a {\em compact Lie algebra}, i.e. any of the following equivalent conditions hold:
 \begin{enumerate}
\item[\ri] $\ggo$ is the Lie algebra of a compact Lie group.
\item[\rii] The Killing form $B$ is negatively semidefinite.
\item[\riii] $\ggo=\overline{\ggo} \oplus \cc$ where  $\cc$ is the center of $\ggo$ and $\overline{\ggo}=C^1(\ggo)$  denotes a semisimple Lie algebra, whose Killing form is negative definite.
\end{enumerate}

The previous paragraph points out that the existence of an ad-invariant metric has a relationship with the algebraic structure of $\ggo$. Bordemann proved in \cite{B} a result which establishes that an ad-invariant metric exists on  $\ggo$  if and only if the adjoint and the coadjoint representation are equivalent. Here we shall give a proof of a weaker result useful for our purposes.

 % The following proposition gives an equivalent condition for Lie algebras admitting an ad-invariant metric. Recall that the coadjoint representation $\ad^{\ast}_u:\ggo^{\ast} \to \ggo^{\ast}$, where $\ggo$ denotes the dual space of $\ggo$, is defined as
 %$$\ad_u^{\ast} \xi = -\xi \circ \ad_u \qquad u \in \ggo, \quad \xi \in \ggo^{\ast}.$$
 
 \begin{prop} Let $\ggo$ be a Lie algebra. Then $\ggo$ admits an ad-invariant metric if and only if there exists a linear isomorphism $T: \ggo \to \ggo^{\ast}$ satisfying $T(x) y = T(y) x$ for all $x, y\in \ggo$ such that \begin{equation}\label{coad}
 \ad^{\ast}_u = T \circ \ad_u \circ T^{-1} \qquad \mbox{ for all }  u \in \ggo.
 \end{equation}
 \end{prop}
 \begin{proof} Assume first that $\ggo$ can be endowed with an ad-invariant metric $\la \,,\, \ra$. The isomorphism induced by the metric  given by $T(x) (y) = \la  x, y \ra$ for all $x, y \in \ggo$ satisfies the requirements. Conversely assume that there exists   a linear isomorphism $T: \ggo \to \ggo^{\ast}$ such that $T(x) (y) = T(y) (x)$ and (\ref{coad}) holds. Define a metric on $\ggo$ by $\la x, y\ra = T(x) y$. The  assumptions show that the bilinear map $\la \,,\, \ra$ is a metric and 
 (\ref{coad}) says that  it is ad-invariant, 
 %for all $u,v, x\in \ggo$ it holds
 %$$\ad^{\ast}_u Tv(x) = -Tv \circ \ad_u (x) = -\la v, [u, x]\ra =( T \circ \ad_u v )(x)= \la [u, v], x\ra,$$
 finishing the proof.
 \end{proof}
 
The so called central descending and ascending  series of a Lie algebra $\ggo$, respectively $\{C^r(\ggo)\}$ and $\{C_r(\ggo)\}$ for all $r\geq 0$, are constitued by the ideals in $\ggo$, that for non negative integers $r$ are given by
$$\begin{array}{rclrcl}
C^0(\ggo) & = & \ggo & C_0(\ggo) & = & 0 \\
C^r(\ggo) & = & [ \ggo, C^{r-1}(\ggo)] \quad & C_r(\ggo) & = & \{ x\in \ggo : [x, \ggo]\in C_{r-1}(\ggo)\}
\end{array}
$$
 Note that $C_1(\ggo)$ is by definition the center of $\ggo$, that will be denoted with $\zz(\ggo)$ or simply $\zz$.
 
Lie algebras for which there exists a positive integer $r$ such that $C^r(\ggo)=0$ are called {\em nilpotent}. It is said {\em k-step nilpotent} if $C^k(\ggo)=0$ but $C^{k-1}(\ggo) \neq 0$.

\begin{lem}\label{lem2} Let $(\ggo, \la \,, \, \ra)$ be a Lie algebra equipped with an ad-invariant metric. Then 

\begin{enumerate}
\item[\ri] If $\hh$ is an ideal then $\hh^{\perp}$ is also an ideal in $\ggo$.

\item[\rii] $C^r(\ggo)^{\perp} = C_r(\ggo)$ for all $r\in \NN_0$.
%If $\jj$ is a one dimensional ideal then it is contained in the center of $\ggo$.
\end{enumerate}
\end{lem}
\begin{proof} The first assertion follows by a direct computation. The second one can be proved by an inductive reasoning. Indeed ii) is clearly true for $r=0$. Assume now that $C^k(\ggo)^{\perp} = C_k(\ggo)$. If $x\in C^{k+1}(\ggo)^{\perp}$  we have $0=\la x, [\ggo, C^k(\ggo)]\ra= \la [x, \ggo], C^k(\ggo)\ra$ and this says that $[x, \ggo] \in C^k(\ggo)^{\perp} = C_k(\ggo)$, therefore $x \in C_{k+1}(\ggo)$. It is not difficult to see that any $x\in C_{k+1}(\ggo)$ belongs to $C^{k+1}(\ggo)^{\perp}$.
%For the general proof see \cite{Me}, we shall prove that ${C^1(\ggo)}^{\perp}=\zz$. If $z\in \zz$. one easily check that $z\in {C^1(\ggo)}^{\perp}$. Assume now $x\in C^1(\ggo)$, thus $0= \la x, [u,v]\ra = \la v,[x,u]\ra$ for all $u,v\in \ggo$. Therefore $[x,u]=0$ for all $u \in \ggo$ and this says that $x\in \zz$. 
\end{proof}

As  corollary of the last lemma on any Lie algebra $\ggo$  admitting an ad-invariant metric the following equality holds 
$$\dim \ggo = \dim \zz + \dim C^1(\ggo).$$
In particular if $\ggo$ is solvable  its center must be non trivial. But this condition is not sufficient for a Lie algebra to admit an ad-invariant metric. To exemplify this assertion let $\nn$ be the 2-step nilpotent Lie algebra $\nn=\zz \oplus \vv$ where $\zz=span\{z_1, z_2, z_3, z_4\}$  and $\vv=span\{v_1, v_2, v_3, v_4\}$ with the non trivial Lie bracket relations
$$[v_1,v_2]= z_3,\quad [v_2, v_3]= z_4, \quad [v_3, v_4]= z_1,\quad [v_1, v_4]=z_2.$$ 
 Then $\zz=C^1(\nn)$ and its dimension coincides with the dimension of $\vv$ but simple computations show that any ad-invariant symmetric bilinear form on $\nn$ is degenerate (see (\ref{four})).

 \section{the structure of 2-step  nilpotent Lie algebras admitting an ad-invariant metric} 
 
 In this section we shall study the algebraic structure of 2-step nilpotent Lie algebras $\nn$ admitting an ad-invariant metric. We shall prove that  the existence of such metric imposes strong conditions on the algebraic structure of $\nn$.
  
 If $\nn$ denotes a 2-step nilpotent Lie algebra the commutator is contained in the center $C^1(\nn)\subset \zz$ and thus we define the  {\em corank} of $\nn$  as the non negative integer $k$ such that  $k:=\dim \zz - \dim C^1(\ggo)$.

 Examples of 2-step nilpotent Lie algebras with an ad-invariant metric can be achieved as in (\ref{cot}) starting with a 2-step nilpotent Lie algebra $\nn$, that is, the cotangent of $\nn$, $\nn^{\ast} \rtimes \nn$, is  also a 2-step nilpotent Lie algebra that admits an ad-invariant metric provided by the canonical hyperbolic one. 
 
 The following proposition shows a orthogonal splitting of any 2-step nilpotent Lie algebra $\nn$ equipped with an ad-invariant metric.

 \begin{prop}\label{p1} Let $(\nn, \la \,, \,\ra )$ be a 2-step nilpotent Lie algebra with an ad-invariant metric. Then $\nn$ is a orthogonal direct sum
 $$ \nn = \tilde{\zz} \times \tilde{\zz}^{\perp}$$
 where $\tilde{\zz}$ is a non degenerate central ideal and $\tilde{\zz}^{\perp}$ is a 2-step nilpotent ideal of corank zero.
 \end{prop}
 \begin{proof} Since $C^1(\nn) = \zz^{\perp} \subset \zz$ any complementary subspace of $C^1(\nn)$ in $\zz$ must be non degenerate. Let $\tilde{\zz}$ denote a fixed complementary subspace, thus
 $$\nn = \tilde{\zz} \oplus \tilde{\zz}^{\perp}$$
 The central ideal $\tilde{\zz}$ is by construction non degenerate and therefore  $\tilde{\zz}^{\perp}$ is a non degenerate 2-step nilpotent ideal (Lemma (\ref{lem2})), whose center  coincides with its commutator $C^1(\tilde{\zz}^{\perp})=C^1(\nn)$; therefore $\tilde{\zz}^{\perp}$ has corank zero.  The ad-invariant metric of $\nn$ induces an ad-invariant metric on $\tilde{\zz}^{\perp} \simeq \nn/\tilde{\zz}$ by restriction. Moreover  since $C^1(\tilde{\zz}^{\perp})=C^1(\nn)$ is a totally isotropic subspace, by applying Lemma (\ref{l1}) to  $\tilde{\zz}^{\perp}$  it is possible to find a totally isotropic complementary subspace $\vv$ such that $\tilde{\zz}^{\perp} = C^1(\nn) \oplus \vv$. Therefore the 2-step nilpotent Lie algebra $\nn$ decomposes as a orthogonal product 
 $$\nn = \tilde{\zz} \times ( C^1(\nn) \oplus \vv )$$ 
 The vector space $\vv$  endowed with the trivial Lie bracket  is isomorphic to $\nn / \zz$.
 \end{proof}

 Among other constructions, 2-step nilpotent Lie algebras with an ad-invariant metric  can be obtained in the following way.  Let $\vv$ denote a vector space  equipped with an inner product $\la \,, \,\ra_{\vv}$  and let  $\rho : \vv \to \sog(\vv)$ denote an injective linear map satisfying 
 \begin{equation}\label{ss}
 \rho(u) v + \rho(v) u = 0 \qquad \mbox{for all }u,v \in \vv
 \end{equation}
 or equivalently $\rho(x)x =0$ for all $x\in \vv$. Let $\nn(\vv, \rho)=\vv^* \oplus \vv$ be the vector space endowed with  the canonical hyperbolic metric $\la \,, \,\ra$ as in (\ref{cot}), that is 
 \begin{equation}\label{hyp}
 \la z_1 + v_1 , z_2 + v_2 \ra = z_1 \cdot v_2 + z_2 \cdot v_1 \qquad \mbox{ for }  z_i \in \vv^*, \, v_i \in \vv, \, i=1,2
 \end{equation}
 where $z \cdot v$ denotes the evaluation. Define a Lie bracket $[\cdot , \cdot ]$ on $\nn(\vv, \rho)$ by 
  $$\begin{array}{rcl}
[\vv^{\ast}, \nn(\vv, \rho)] & = & 0  \quad \mbox{and }\quad [\vv, \vv ] \subset \vv^{\ast} \quad \mbox{ with } \\
\la [u,v], w\ra & =& \la \rho(w) u, v \ra_{\vv}\quad \mbox{ for all } \quad u,v , w \in \vv.
\end{array}
$$
Since $\rho(u)\in \sog(\vv)$ for all $v\in \vv$ then $[\cdot ,\cdot ]$ is skew-symmetric and (\ref{ss}) implies that $\ad_x$ are skew-symmetric with respect to $\la \,, \,\ra$. Thus $\nn(\vv, \rho)$  is a 2-step nilpotent  Lie algebra  with an ad-invariant metric of corank zero. In fact, since  $\rho$ is injective we have $C^1(\nn(\vv, \rho)) = \vv^{\ast} = \zz$. We call the  2-step nilpotent Lie algebra $\nn(\vv, \rho)$ above the {\em modified cotangent} of $\vv$.  

 Finally if we provide $\RR^m$ with a metric and take $\RR^m \oplus \nn(\vv, \rho)$ the orthogonal direct product of the Lie algebras with ad-invariant metrics, then we obtained a 2-step nilpotent Lie algebra of corank $m$ provided with an ad-invariant metric. 
 
 We say that two Lie algebras $(\ggo_1, \la \,, \,\ra_1)$ and $(\ggo_2, \la \,, \,\ra_2)$ are {\em isomorphic isometric} if there exists an isomorphism of Lie algebras $\varphi:\ggo_1 \to \ggo_2$ which is also an isometry.
 
 The following theorem shows that any 2-step nilpotent Lie algebra endowed with an ad-invariant metric is isomorphic isometric to a direct product of an abelian Lie algebra  and a modified cotangent of the vector space $\vv$, where $\vv= C^1(\nn)^{\ast}$.

  \begin{thm}\label{t1} Let $({\nn}, \la \,, \,\ra)$ denote a 2-step nilpotent Lie algebra  of rank $m$ endowed with an ad-invariant metric. Then $({\nn}, \la \,, \,\ra)$ is isomorphic isometric to a orthogonal direct product of the Lie algebras $\RR^m$ and the modified cotangent of $C^1(\nn)^*$. 
   \end{thm}
  \begin{proof} By Proposition (\ref{p1}) the Lie algebra $\nn$ decomposes into a orthogonal direct sum $\nn=\tilde{\zz} \oplus \tilde{\zz}^{\perp}$ where $\tilde{\zz}$ is a non degenerate abelian ideal (and hence isomorphic isometric to $\RR^m$ for some $m$) and  $\tilde{\zz}^{\perp}$ is a non degenerate  2 -step nilpotent ideal of $\nn$. Since  $C^1(\nn) =C^1(\tilde{\zz}^{\perp})$ is a totally isotropic ideal of $\tilde{\zz}^{\perp}$,  there exists a totally isotropic complementary subspace $\vv$ such that $\tilde{\zz}^{\perp} = C^1(\nn) \oplus \vv$ and clearly $C^1(\nn)$ is isomorphic to $\vv^{\ast}$ as  vector spaces (via the metric on $\nn$). Equipp $\vv$ with a inner product $\la \,,  \,\ra_{\vv}$ and define a linear map $\rho:\vv \to \sog(\vv)$ by
  $$\la \rho(u) v, w \ra_{\vv} = \la [ v, w], u\ra \qquad \mbox{ for all } u,v , w \in \vv$$
  Indeed $\rho(\vv) \subset \sog(\vv, \la \,, \,\ra_{\vv})$ since the Lie bracket on $\nn$ is skew-symmetric and the ad-invariance property of the metric implies that  $\rho$ satisfies (\ref{ss}). Let $\nn(\vv, \rho)$ denote the modified cotangent of $\vv$ defined as above,  attached with the Lie bracket  $[ \cdot, \cdot]'$ and hyperbolic metric $\la \,, \, \ra'$. 
   
  Assume that $\{z_1, \hdots, z_n\}$ is a basis of $C^1(\nn)$ and let $\{v_1, \hdots, v_n\}$ be a basis of $\vv$ such that $\la z_i, v_j\ra =\delta_{ij}$. Take $\{\phi_1, \hdots, \phi_n\}$ the dual basis of $\vv$ in $\vv^*$. Let $\varphi: \zz^{\perp} \to \nn(\vv, \rho)$ be the isometry given by $\varphi(z_i)=\phi_i$ and $\varphi(v_i)=v_i$. Then $\varphi$ is also a automorphism of Lie algebras. In fact
  $$\la [\varphi v, \varphi w]',  \varphi u\ra'= \la \rho(u) v, w \ra_{\vv} = \la [v, w], u\ra= \la \varphi [v, w]', \varphi u\ra'.$$
   Finally the central ideal  $\tilde{\zz}$  is isomorphic to $\RR^m$ where $m=\dim \tilde{\zz}$ and they are also isometric by transporting  the metric. 
   \end{proof}

   In view of the previous theorem we shall denote a 2-step nilpotent Lie algebra endowed with an ad-invariant metric as $\nn(\vv, \rho, \la \,, \,\ra_{\vv})$.
  
   \begin{rem} Another version of the previous theorem can be done by replacing the totally isotropic vector space $\vv$ by a non degenerate  vector space $\vv'$ such that $C^1(\nn)\oplus \vv=C^1(\nn)\oplus \vv'$ by applying Lemma (\ref{l1}). Furthermore $\vv'$ can be choosen in such way that the restriction of the metric on $\tilde{\zz}^{\perp}$ defines a inner product. In this way, one gets an injective homomorphism satisfying (\ref{ss}) $\rho:\vv' \to \sog(\vv', \la \,, \,\ra'_{\vv'})$ where $\la \,, \,\ra'_{\vv'}$ is the restriction of the metric of $\nn$.
    \end{rem}
 
Assume $\nn$ is a 2-step nilpotent Lie algebra of corank zero admitting an ad-invariant metric $\la \,, \,\ra$ and let $(\vv, \la \,, \,\ra_{\vv})$ be a  isotropic vector subspace of $\nn$ complementary to the commutator  $C^1(\nn)$ and suppose $\vv$ is equipped with a inner product $\la \,, \,\ra_{\vv}$. Let $\rho: \vv \to \sog(\vv)$ denote a injective map satisfying $\rho(u)v +  \rho(v) u = 0$ for all $u,v \in \vv$. Suppose $\{v_i\}$ is a orthonormal basis of $\vv$ for i=$1, \hdots, p$ with $z_1, \hdots z_p$ a dual basis of $\zz$ such that $\la z_i, v_j\ra=\delta_{ij}$. As in (\ref{t1}) 
 $\la \rho(v_i) v_j, v_k\ra_{\vv} = \la [v_j,v_k], v_i \ra \qquad \text{for all }u,v,w \in \vv.$ Therefore 
denoting   the image by $\rho$ of $v_i$ by $A^i = \rho(v_i)$   and by $a_{jk}^i$ the $jk$-entry of $A^i$, we have
$$a^i_{jk}= \la A^i v_j, v_k\ra_{\vv} = \la [v_j,v_k], v_i \ra$$
Thus   the coefficients  $a_{jk}^i$ satisfy

\begin{equation}\label{coef}
\begin{array}{llll}
a^i_{kj} & = & - a^i_{jk} & \mbox{ since $A_i\in \sog(\vv)$ for all i;} \\
a^i_{kj} & = & - a^j_{ki} & \mbox{ for all i,j,k, which follows from (\ref{ss}).}
\end{array}
\end{equation}

 Both conditions  imply that there is no 2-step nilpotent Lie algebra of corank zero admitting an ad-invariant metric  with one or two-dimensional commutator. In the three dimensional case  the maps  $\rho(v_i)\in so(\vv)$ can be represented  by the following matrices
$$
A^1=\left( \begin{matrix}
0 & 0 & 0  \\
0 & 0 & a \\
0 & -a & 0 
\end{matrix} \right)  \qquad 
A^2= \left( \begin{matrix}
0 & 0 & -a \\
0 & 0 & 0 \\
a & 0 & 0  
\end{matrix} \right)\qquad 
A^3= \left( \begin{matrix}
0 & a & 0 \\
-a & 0 & 0 \\
0 & 0 & 0 
\end{matrix} \right)
$$
which are linearly independent for any non zero real number $a$. Thus  the set of 2-step nilpotent Lie algebras of corank zero and three dimensional commutator is parametrized by $\RR -\{0\}$. The following result explains the existence  in any  dimension.

\begin{prop} \label{four} If the dimension of the vector space $\vv$ is three or greather than four, there exists a modified cotangent $\nn(\vv, \rho)$. 
\end{prop}
\begin{proof} The explanations above show that there is no 2-step Lie algebra of corank zero admitting an ad-invariant metric and whose commutator has dimension one or two. If the dimension of $\vv$  is three, we showed above the form of the matrices $A^1, A^2, A^3$ satisfying (\ref{coef}) to construct  a 2-step nilpotent Lie algebra as a modified cotangent.  Let $\vv$ be  a four dimensional vector space and fix a inner product $\la \, , \, \ra$ on $\vv$. Assume that $\{v_1, v_2, v_3, v_4\}$ is a orthonormal basis for $\la \,, \,\ra_{\vv}$. The set of linear transformations in $\sog(\vv)$ satisfying (\ref{coef})  can be represented by the following matrices
$$
A^1=\left( \begin{matrix}
0 & 0 & 0 & 0 \\
0 & 0 & a & b \\
0 & -a & 0 & c \\
0 & -b & -c & 0
\end{matrix} \right)  \qquad 
A^2= \left( \begin{matrix}
0 & 0 & -a & -b \\
0 & 0 & 0 & 0 \\
a & 0 & 0 & d \\
b & 0  & -d & 0
\end{matrix} \right)$$
$$ 
A^3= \left( \begin{matrix}
0 & a & 0 & -c\\
-a & 0 & 0 & -d \\
0 & 0 & 0 & 0 \\
c & d & 0 & 0
\end{matrix} \right)  \qquad
A^4= \left( \begin{matrix}
0 & b & c & 0 \\
-b & 0 & d & 0 \\
-c & -d & 0 & 0 \\
0 & 0 & 0 & 0
\end{matrix} \right). 
$$  
It is not difficult to prove that this set cannot be linearly independent for any choice of the scalars $a,b,c,d$. Hence any linear map from $\vv$ to $\sog(\vv)$ satisfying (\ref{ss}) cannot be injective.

In what follows we shall show the existence of a set of  n linearly indepent matrices $A^i \in \sog(\vv)$ in  any dimension $n \geq 5$ satisfying (\ref{coef}).  We shall show a contruction of this set. Indeed this is not unique. 

Assume $n\equiv 0$ (mod 3).  Let $A^i$ be a block matrix consisting of 3$\times$3 matrices, all of them zero with exception of one in the principal diagonal of $A^i$:
 $$
A^i= {\left( \begin{matrix}
0 &  & 0 & & 0\\
0 &\ddots  & 0 & & 0\\
0 &   & r^i &  & 0\\
0 &  & 0 & \ddots & 0 \\
0 &   & 0 &   & 0
\end{matrix} \right)\qquad }
$$ 
 
 with  $r^i$, the 3 $\times$ 3 block in the k+1 place of the diagonal, $k=\frac{i-s}3$ s=0,1,2, given by
$$\begin{array}{lll}
r^{3 k +1}={\left( \begin{matrix}
0 & 0 & 0 \\
0 & 0 & -1\\
0 & 1 & 0 
\end{matrix} \right)\qquad } & 
r^{3 k +2}={\left( \begin{matrix}
0 & 0 & 1 \\
0 & 0 & 0\\
-1 & 0 & 0 
\end{matrix} \right)\qquad } & 
r^{3 k +3}={\left( \begin{matrix}
0 & -1 & 0 \\
1 & 0 & 0\\
0 & 0 & 0 
\end{matrix} \right)} \\
{\mbox{if }\,i= 3 k +1} & {\mbox{if }\, i=3k+2} & {\mbox{if }\, i=3k+3}
\end{array}
$$

If $n\equiv 1$ (mod 3) or $n \equiv 2$ (mod 3) we define matrices $A^i$ as above for all i=1, $\hdots, n-3$, that is $A^i$  consists on zero $3\times 3$ block matrices  except for the $r$-block on the principal diagonal being $i=3(r-1)+s$, and s=1,2,3 as above. For $n-2, n-1, n$ take $A^{n-2}, A^{n-1}$ and $A^n$ respectively as the matrices consisting on zero blocks with exception of the last $5\times 5$ block on the diagonal that we choice as follows:
$$\begin{array}{lll}
{\left( \begin{matrix}
0 & -1 & 0 & 0 & 0\\
1 & 0 & 0 & 0 & 0\\
0 & 0 & 0 & 0 & 0 \\
0 & 0 & 0 & 0 & -1 \\
0 & 0 & 0 & 1 & 0 
\end{matrix} \right)\qquad } & 
{\left( \begin{matrix}
0 & 0 & 0 & 0 & 0\\
0 & 0 & 0 & 0 & 0\\
0 & 0 & 0 & 0 & 1 \\
0 & 0 & 0 & 0 & 0 \\
0 & 0 & -1 & 0 & 0 
\end{matrix} \right)\qquad } & 
{\left( \begin{matrix}
0 & 0 & 0 & 0 & 0\\
0 & 0 & 0 & 0 & 0\\
0 & 0 & 0 & 1 & 0 \\
0 & 0 & -1 & 0 & 0 \\
0 & 0 & 0 & 0 & 0 
\end{matrix} \right)} \\
{\mbox{if }\,i= n-2} & {\mbox{if }\, i=n-1} & {\mbox{if }\, i=n}
\end{array}
$$ 

 Since $A^i$ were constructed as block matrices, canonical computations show that the set $\{A^i\}$ is linearly independent.
\end{proof}

Next we shall study the set of injective homomorphisms $\rho : \vv \to \sog(\vv)$ satisfying $ \rho(u) v + \rho(v) u = 0$. Indeed this is a subset of the $n^2(n-1)/2$-dimensional vector space $Hom(\vv, \sog(\vv))$,  where $n=\dim \vv$.

\begin{prop} Let $(\vv, \la \,, \,\ra_{\vv})$ denote a vector space equipped with a inner product $\la \,, \,\ra_{\vv}$ and let $\mathcal R$ denote the set
 $$\mathcal R= \{ \rho:\vv \to \sog(\vv) \mbox{ injective linear map  satisfying } \rho(u) v + \rho(v) u = 0\}$$ Then  $\mathcal R$ is closed in $Hom(\vv, \sog(\vv))$, moreover $\mathcal R$ is a subset of the open set of homomorphism of rank n.
\end{prop}
\begin{proof} The set $\mathcal O$ consisting of $t\in Hom(\vv,\sog(\vv))$ of rank n is  open  and $\mathcal R \subset \mathcal O$ since $\rho$ is injective. Fix a vector $u\in \vv$ and consider $\chi_u:Hom(\vv, \sog(\vv)) \to Hom(\vv)$ the map given by 
$$\chi_u(\rho)(v)= \rho(u) v + \rho(v) u$$
Clearly $\chi_u(\rho)$ is a linear map over $\vv$ and $\chi_u$ is a linear  morphism  with kernel $\kk(u)=\{\rho:\vv \to \sog(\vv) / \rho(u) v + \rho(v) u=0 \mbox{ for all } v \in \vv\}$ which is a closed subspace of $Hom(\vv, \sog(\vv))$. Since
$$\mathcal R =\cap_{u \in \vv} \kk(u)$$
then  $\mathcal R$ is closed in $Hom(\vv, \sog(\vv))$. 
\end{proof}

\subsection{Isomorphism classes} In this section we investigate the isomorphism and isometric isomorphism classes of 2-step nilpotent Lie algebras with ad-invariant metrics.

Let $\nn(\vv, \la \,, \,\ra, \rho)$ and $\nn'(\vv', \la \,, \,\ra', \rho')$ denote 2-step nilpotent Lie algebras endowed with ad-invariant metrics (\ref{t1}). Since any isomorphism of Lie algebras preserves the center we deduce that if $\vv$ is not isomorphic to $\vv'$ then $\nn(\vv, \la \,, \,\ra, \rho)$ cannot be isomorphic to $\nn'(\vv', \la \,, \,\ra', \rho')$. 

We shall first fix the vector space $\vv$ and we shall study the isomorphism classes. As we prove below the isomorphism class is independent of the inner product on $\vv$.

\begin{prop} \label{p3} Let $\la \,, \,\ra$ and $\la \,, \,\ra'$ be inner products on a vector space $\vv$ and assume the injective linear map $\rho:\vv\to \End(\vv)$ satisfying (\ref{ss}) is skew-symmetric with respect to both inner products. Then the corresponding 2-step nilpotent Lie algebras $\nn$ and $\nn'$ constructed up $(\la \,, \,\ra, \rho)$ and $(\la \,, \,\ra', \rho)$ respectively are isomorphic.
\end{prop}
\begin{proof} Let $P:\vv \to \vv$ denote the positive definite transformation such that 
\begin{equation}\label{ps}\la x,y \ra' = \la Px, y \ra = \la P^{1/2} x, P^{1/2} y \ra
\end{equation}
holds. Then  for all $u,v,w \in \vv$:
$$\la P \rho(u) v, w \ra = \la \rho(u)v , w \ra' = - \la v, \rho(u) w\ra' = -\la P v, \rho(u) w \ra.$$
and this implies that $P \rho(u) = \rho(u) P$ for all $u\in \vv$. Thus the only symmetric square root of $P$ denoted by $P^{1/2}$ also satisfy $P^{1/2} \rho(u) = \rho(u) P^{1/2}$. We then have that $(P^{\ast}, P^{1/2}):\nn'=\vv^{\ast} \oplus \vv \to \nn=\vv^{\ast} \oplus \vv$ is an isomorphism of Lie algebras. Indeed if $u,v,w\in \vv$ then 
$$\begin{array}{rcl}
\la [P^{1/2} u, P^{1/2}v]_{\nn}, w \ra & = & \la \rho(w) P^{1/2} u, P^{1/2} v\ra = \la P^{1/2} \rho(w) P^{1/2} u,  v\ra \\
& = & \la P \rho(w)  u,  v\ra =  -\la P\rho(u) w, v\ra' \\
& = & -\la \rho(u) Pw, v \ra'= \la \rho(Pw) u, v\ra \\
& = & \la [u,v]_{\nn'}, w \ra' = \la P^{\ast}[u,v]_{\nn'}, w \ra
\end{array}
$$ where ${P}^{\ast}:\vv^{\ast}_{\la \,, \,\ra} \to \vv^{\ast}_{\la \,, \,\ra'}$ is  the dual of $P:\vv_{\la \,, \,\ra'} \to \vv_{\la \,, \,\ra}$.
\end{proof}

Let  $\vv$ denote  vector space endowed with inner products $\la \,,\, \ra_{\vv}$, $\la, \ra'_{\vv}$ and assume $[\cdot , \cdot]'$ is  a Lie bracket on $(\nn' = \vv^{\ast} \oplus \vv, \la \,, \,\ra'_{\nn})$, with respect to which  $\la \,, \,\ra'_{\nn}$ is ad-invariant. It  is always possible to define a Lie bracket on $(\nn = \vv^{\ast} \oplus \vv, \la \,, \,\ra_{\nn})$  for which $(\nn, \la \,, \,\ra)$ and $(\nn', \la \,, \,\ra')$ are isomorphic isometric. This can be done by transporting the Lie bracket $[\cdot , \cdot ]'$ via the inner product. In fact
let $P:\vv \to \vv$ denote the positive definite transformation as in (\ref{ps}) such that
\begin{equation}\la x,y \ra' = \la Px, y \ra = \la P^{1/2} x, P^{1/2} y \ra
\end{equation}
where $P^{1/2}$ denotes the only symmetric square root of $P$. Transport the Lie bracket $[\cdot , \cdot ]'$ via $P^{1/2}$, that is, set $[x,y]= {P^{-1/2}}^{\ast} [P^{-1/2}x, P^{-1/2}y]'$, where ${P^{-1/2}}^{\ast}:\vv^{\ast}_{\la \,, \,\ra'} \to \vv^{\ast}_{\la \,, \,\ra}$. Then the map $({P^{-1/2}}^{\ast}, P^{1/2}):(\nn, [ \cdot, \cdot]', \la \,, \,\ra'_{\nn}) \to (\nn, [\cdot ,\cdot ], \la \,, \,\ra_{\nn})$ is an  isomorphism of the Lie algebras $\nn$ and $\nn'$ since ${P^{-1/2}}^{\ast}[x,y]'= {P^{-1/2}}^{\ast} [ P^{-1/2} P^{1/2}x, P^{-1/2}P^{1/2}y]'= [P^{1/2}x, P^{1/2}y]$. Furthermore it is an isometry.

In view of the previous results we shall fix the vector space $\vv$ with the inner product  and  we  denote the 2-step nilpotent Lie algebra of corank zero of Theorem (\ref{t1}) constructed with the data $(\vv , \la \,, \,\ra_{\vv}, \rho)$ just as $\nn(\vv,\rho)$. 

Our goal now is to describe the isomorphisms between 2-step nilpotent Lie algebras equipped with ad-invariant metrics. 

\begin{thm}\label{t2} Let $\nn$ and $\nn'$ denote the 2-step nilpotent Lie algebras constructed from a vector space $(\vv,\la \, , \, \ra_{\vv})$ with respective homomorphism $\rho$ and $\rho'$. Let  $\la \,,\,\ra$ denote the  ad-invariant metric on $\nn$. Then $\nn$ is isomorphic to $\nn'$ if and only if {\rm corank} $\nn$ = {\rm corank} $\nn'$ and there exist non singular linear transformations $A,B \in \Gl(\vv)$ such that
\begin{equation}\label{iso}
\rho(A v) = B^t  \rho'(v) B \qquad \quad \mbox{ for all } v \in \vv,
\end{equation}
where $B^t$ denotes the transpose for the inner product $\la \,, \,\ra_{\vv}$ of  $B$ given by $\la B^t u, v\ra_{\vv}= \la u, Bv \ra$.

If $\tilde{\zz}$ is a complementary subspace of $\vv^*$ in $\zz$ and $\la \,, \,\ra_{\tilde{\zz}}$ denotes the restriction of the metric of $\nn$ to  $\tilde{\zz}$ then any isometric isomorphism $\phi:\nn \to \nn'$ can be represented by a matrix as follows:
\begin{equation}\label{oiso}
\left( \begin{matrix} C & 0 & D \\ E & {A^t}^{-1} & D'\\ 0 & 0 & A^{} \end{matrix} \right)
\qquad {\begin{array}{l} \mbox{ with } A\in \Gl(\vv), C \mbox{ an isometry for }\la  \,, \, \ra_{\tilde{\zz}} \mbox { and }\\
C^t D + E^tB=0 \quad  D^t D + {D'}^t B + B^t D' =0 \end{array}}
\end{equation}
\end{thm}
\begin{proof} Using Theorem (\ref{t1}) assume $\nn\simeq \RR^m \times \nn(\vv, \rho)$ and $\nn' \simeq \RR^{m'} \times \nn(\vv, \rho)$. Since the  isomorphism  $\phi:\nn \to \nn'$ leaves the center and the commutator invariant, $\phi$ induces a isomorphism  $\tilde{\phi}:\zz(\nn)/C^1(\nn)\to \zz(\nn')/C^1(\nn')$, therefore $corank \, \nn = corank \,\nn'$.  Moreover there exist linear operators  $A, B\in \Gl(\vv)$ such that $A^{\ast} [u,v]=[Bu, B v]$ where $A^{\ast} : \vv^{\ast} \to \vv^{\ast}$ is the dual of a non singular linear map  $A:\vv \to \vv$. In fact if $\tilde{\zz}$ denotes a  complementary subspace of the commutator in the center, for every $v\in \vv$ assume $\phi v= Dv +D'v+ Bv$, where $D:\vv \to \tilde{\zz}$, $D':\vv \to \vv^*$ and $B: \vv \to \vv$. Therefore  $\phi [u,v] = [\phi u, \phi v]'$ if and only if
$$ \la A^{\ast}[u,v], w \ra = \la [Bu, Bv]', w \ra$$ 
if and only if $$\la \rho(A w)u,v \ra_{\vv} = \la \rho'(w) Bu, Bv  \ra_{\vv}= \la B^t\rho'(w) Bu, Bv  \ra_{\vv}$$ 
and so we get condition (\ref{iso}). The isomorphism $\phi$ can be represented by  a matrix of the form
$$
\left( \begin{matrix} C & 0  & D \\ E & A^{\ast} & D' \\ 0 & 0 & B \end{matrix} \right)
\qquad \mbox{ with } A,  C \in \Gl(\vv) \quad \mbox{and } %\quad C^t B + B^t C + C^t C = I
$$
relative to a basis adapted to the decomposition $\tilde{\zz} \oplus \vv^* \oplus \vv$.
Conversely  assume  $C:\tilde{\zz} \to \tilde{\zz}$ is a linear isomorphism and take $A,B \in Gl(\vv)$ satisfying (\ref{iso}). Define $\phi:\nn \to \nn'$  by $\phi(v) = Bv$ for $v\in \vv$,  $\phi(z)= Cz$ for $z\in \tilde{\zz}$ and $\phi(x) = A^{\ast} x$ for $x \in \vv^{\ast}=C^1(\nn)$, being  $A^{\ast} : \vv^{\ast} \to \vv^{\ast}$ the dual of $A$. Then it is not hard to prove that $\phi$ is an isomorphism of Lie algebras. 

If we assume  that a isomorphism of Lie algebras $\phi$ is orthogonal, the matrices $A,B, C$, $D, D', E$ must satisfy:
$$ C^t D + E^t B =0\qquad \qquad D^t D + {D'}^t B + B^t D'=0  $$
where $C$ must be orthogonal relative to $\la \,, \,\ra_{\tilde{\zz}}$ and $AB=I$, therefore  $\phi$ can be represented by a matrix of the form  (\ref{oiso}).
\end{proof}

Notice that the matrix in (\ref{oiso}) reduces to a simple form if $\nn$ has corank zero. In this case we get $C=0$, $D=0$, $E=0$ and  we must ask ${D'}^t A^{-1} + {A^{-1}}^t D' =0$.

\begin{cor} There exists only one six-dimensional 2-step nilpotent Lie algebra of corank zero admitting an ad-invariant metric. 
\end{cor}

\begin{exam} In this example we show that a 2-step nilpotent Lie algebra $\nn$ with ad-invariant metric $\la \, , \, \ra$  can be achieved as a modified cotangent in different ways. 
 
 Let $\RR^3$ equipped with its canonical inner product $( \,,\, )$. Let $T: \RR^3 \to \sog(3)$ be the linear transformation defined by
 $$
 T(x_1, x_2, x_3) = \left( \begin{matrix}
0 & -x_3 & x_2  \\
x_3 & 0 & -x_1 \\
-x_2 & x_1 & 0 
\end{matrix} \right).$$ 
Then $T$ applies $\RR^3$ onto $\sog(3)$ and it holds $T(u) v + T(v) u=0$ for all $u,v \in \RR^3$. But $T$ cannot be a representation of $\RR^3$ (considered with the trivial Lie bracket) since $AB \neq BA$ for any pair of matrices $A,B\in \sog(3)$. Let $\nn(\RR^3, T)$ be the 2-step nilpotent Lie algebra with ad-invariant metric as in (\ref{t1}) constructed with this data.

 Now  $\sog(3)$ endowed with the inner product given by minus its Killing form $B( U, V ) = -  tr(\ad_U \ad_V)$ for $U,V \in \sog(3)$.  The adjoint representation in $\sog(3)$ satisfies (\ref{ss}). Let $\nn(\sog(3), \ad^{\ggo})$ denote the 2-step nilpotent Lie algebra with ad-invariant metric as in the modified cotangent model (\ref{t1}).
 One easily  verifies that $(u, v) = \frac 12 B( Tu, Tv)$ for all $u,v\in \RR^3$. And the 2-step nilpotent Lie algebras $\nn(\RR^3, T, (\,, \,) )$ and $\nn(\sog(3), \ad^{\ggo}, -B)$ are isometric isomorphic. In fact $\frac{T}{\sqrt2}$ is an isometry and $T(u) = T^{-1} \ad_{Tu} T$ for all $u\in \RR^3$, implying the isometric isomorphism.
\end{exam}

\section{Examples and Applications}

In this section we shall study the 2-step nilpotent Lie algebras of corank zero which arise from representations of semisimple compact Lie algebras, that is, we shall consider the 2-step nilpotent Lie algebras $\nn(\ggo, \pi)$ of zero corank as in Theorem (\ref{t1}) where
\begin{enumerate}
 \item[\ri] $\ggo$ is a semisimple compact Lie algebra endowed with an inner product $\la \,, \, \ra_{\ggo}$.
 
 \item[\rii] $\pi:\ggo \to \sog(\ggo)$ is a faithful representation satisfying $\pi(x) y + \pi(y) x=0$ for all $x,y    \in \ggo$.
 \end{enumerate}
 
 Explicitely, let $(\ggo, [\cdot , \cdot]_{\ggo})$ be a semisimple compact semisimple Lie algebra equipped with a inner product $\la \, , \, \ra_{\ggo}$. Let $\pi: \ggo \to \End(\ggo)$ denote a faithful representation satisfying (\ref{ss}). The set of this kind of representations is not empty. In fact, if $B$ denotes the Killing form on $\ggo$, the adjoint representation $\ad^{\ggo}: \ggo \to \End(\ggo)$ is faithful, the maps $\ad_x^{\ggo}=[ x,  \cdot ]_{\ggo}$ are skew symmetric with respect to $B$ and they satisfy (\ref{ss}). Hence the modified cotangent of $\ggo$, $\nn= \ggo^{\ast} \oplus \ggo$ endowed with its canonical hyperbolic metric, inherits Lie bracket $[\cdot , \cdot ]$ for which $\nn$ becomes  a 2-step nilpotent Lie algebra  that admits an ad-invariant metric. Roughly speaking take $\nn =\ggo_1 \oplus \ggo_2$ where $\ggo_1 =\ggo = \ggo_2$ and $\ggo_1 = \zz(\nn)$ and the Lie bracket is given as follows: for $x,y \in \ggo_2$ take its Lie bracket on $\ggo$, $[x,y]_{\ggo}$, and choose the copy of this in $\ggo_1$.  Moreover  the family $\ad_x^{\ggo}$ for $x\in \ggo$ is a subalgebra in $\sog(\ggo)$.

 \vskip 6pt
 
We fix a compact Lie algebra $\ggo$ and we study the isomorphism classes of 2-step nilpotent Lie algebras with ad-invariant metrics $\nn(\ggo, \pi)$ which can be constructed as modified cotangent by using different representations of $\ggo$. If  $[ \cdot, \cdot]_{\ggo}$ denotes the Lie bracket on $\ggo$ the condition of $\pi$ being representation says 
\begin{equation}\label{rep}
\pi([x,y]_{\ggo})= \pi(x) \pi(y) - \pi(y) \pi(x).
\end{equation} 

Let $[ \cdot , \cdot]_{\pi}:\ggo \times \ggo \to \ggo$ be the bracket defined by $[ x, y ]_{\pi} = \pi(x) y$. Since $\pi(x) y +\pi(y) x =0$ for all $x,y\in \ggo$, the bracket $[\cdot , \cdot ]_{\pi}$ is bilinear and skew-symmetric and it satisfies the Jacobi identity if and only if
$$\pi(\pi(x) y ) = \pi(x) \pi(y) - \pi(y) \pi(x).$$
But if this is the case, compairing (\ref{rep}) and the previous equality we get  that $[\cdot ,\cdot ]_{\ggo}= [ \cdot, \cdot]_{\pi}$ since $\pi$ is injective .

\begin{prop} Let $\ggo$ be a compact Lie algebra endowed with a inner product and assume $\pi:\ggo \to \sog(\ggo)$ is a faithful representation satisfying (\ref{ss}) and the following relation 
\begin{equation}\label{pi}
\pi(\pi(x) y ) = \pi(x) \pi(y) - \pi(y) \pi(x) \qquad \mbox{ for all } x,y \in \ggo,
\end{equation}
then $\pi$ is the adjoint representation. 
\end{prop}

\begin{rem} The simply connected Lie group $N$ with Lie algebra $\nn$ arising from the semisimple compact Lie algebra $\ggo$ is an example of a naturally reductive  nilmanifold (\cite{La}), whenever $N$ is equipped with the Riemannian metric arising from left translations from $\nn$ as follows: $\nn=\ggo_1 \oplus \ggo_2$ is equipped with the inner product obtained  by taking minus the Killing form on each sumand and making the sum orthogonal.
\end{rem}

According to Theorem (\ref{t2}) the two representations $\pi, \pi':\ggo \to \sog(\ggo)$ give rise to isomorphic nilpotent Lie algebras $\nn(\ggo, \pi)$ $\nn(\ggo, \pi')$ if and only if there exists $A,C\in \Gl(\ggo)$ such that 
$$\pi(A x) = C^t  \pi'(x) C \qquad \quad \mbox{ for all } x \in \vv$$
where $C^t$ denotes the transpose with respect to the inner product. If we choose $C \in O(\ggo)$ then $C^t=C^{-1}$ and (\ref{iso}) implies that  $A= \pi^{-1}\circ \Ad(C)\circ \pi'$ and since $\pi: \ggo \to \pi(\ggo)$, $\Ad(C): \End(\ggo) \to \End(\ggo)$ and $\pi': \ggo \to \pi'(\ggo)$ are Lie algebra isomorphisms then $A: \ggo \to \ggo$ is an automorphism of $\ggo$. Furthermore if $\ggo$ is semisimple $Aut(\ggo) \subset O(\ggo)$.

\begin{rem} Notice that the isometric isomorphisms where $A\in O(\ggo)$ gives rise to a subgroup of the group of isometric isomorphisms. Compare this with \cite{La}, where the author studied 2-step nilpotent Lie groups equipped with a Riemannian metric that can be constructed from representations of compact Lie algebras. 
\end{rem}

\subsection{Lie bialgebras} Next we prove the existence  of Lie bialgebras structures on some  2-step nilpotent constructed  from a compact semisimple Lie algebra $\ggo$.

Let $\ggo$ denote a Lie algebra and $\delta: \ggo \to \Lambda^2\ggo$ is a 1-cocycle with respect to the adjoint representation such that $\delta^*:\Lambda^2\ggo^* \to \ggo^*$ induces a Lie algebra structure on $\ggo^*$, then $(\ggo, \delta)$ is called a {\em Lie bialgebra} on $\ggo$ (see \cite{D1}). If $G$ is a simply connected Lie group whose Lie algebra is $\ggo$, there is a one-to-one correspondence between multlplicative Poisson structures on $G$ and Lie bialgebra structures on $\ggo$ (see \cite{D2}).

A linear operator $r \in \End(\ggo)$
 is called a {\em classical r-matrix} if the $\ggo$-valued skew symmetric
 bilinear form on $\ggo$ given by
 $$
  [x, y]_r = [r x,y] + [x,ry]
 $$
 is a Lie bracket, that is, it satisfies the Jacobi identity. The vector space underlying $\ggo$ equipped with the Lie bracket $[\cdot , \cdot ]_r$ will be denoted by $\ggo_r$. The next result shows a construction of  Lie bialgebra structures (see \cite{AGMM} \cite{Sts}).

 \begin{thm} Let $\ggo$ be a Lie algebra equipped with an ad-invariant metric $\la \,, \, \ra$ and a classical r-matrix $r$ which is skew-symmetric relative to $\la \, , \, \ra$, then $r$ induces a Lie bialgebra structure on $\ggo$.
 \end{thm}
 
 We shall prove the existence of Lie bialgebra structures making use of the previous result.  
 Let $\ggo$ denote a compact semisimple Lie algebra with its Killing form $B$ and assume $r \in \End \ggo$ is skew symmetric for a classical r-matrix $r$.  Let $(\nn=\ggo^{\ast} \oplus \ggo, \la \,, \,\ra)$ denote the 2-step nilpotent Lie algebra with ad-invariant metric $\la \,, \,\ra$ constructed as a modified cotangent of $\ggo$ up $(\ggo, -B)$ with the adjoint representation (\ref{t1}). Then $r^{\ast}\times r\in \End \nn$ is a classical r-matrix, where $r^{\ast}:\ggo^{\ast} \to \ggo^{\ast}$ is the dual of $r$. In fact let $[\cdot , \cdot ]$ denote the Lie bracket on $\nn$, thus for $z_i+v_i \in \ggo^{\ast} \oplus \ggo = \nn$ we have 
 $$
  [z_1+v_1,z_2+v_2]_{r^{\ast}\times r} = [r v_1,v_2] + [v_1,rv_2]
 $$
 and therefore
 $$[[z_1+v_1,z_2+v_2]_{r^{\ast}\times r}, z_3+v_3]_{r^{\ast}\times r}= 0,$$ 
 since $[r v_1,v_2] + [v_1,rv_2]$ belongs to the center of $\nn$ which is invariant by $r^{\ast}\times r$. This implies the Jacobi identity for $[\cdot , \cdot]_{r^{\ast} \times r}$. It is not difficult to prove that $r^{\ast} \times r$ is skew symmetric relative to the ad-invariant metric $\la \,, \,\ra$ on $\nn=\ggo^{\ast} \oplus \ggo$. This explanation is summarized in  the following result.
 
 \begin{thm} If $\ggo$ is a semisimple compact Lie algebra with a classical r-matrix, then the 2-step nilpotent Lie algebra  $\nn(\ggo, \ad^{\ggo})$ constructed as a modified cotangent of $\ggo$ admits a classical r-matrix. Furthermore if $r$ is skew symmetric relative to the Killing form on $\ggo$ then $\nn(\ggo, \ad^{\ggo})$ can be provided with a Lie bialgebra structure.
 \end{thm}
 
Examples of $r$-matrices are given by complex structures. In fact if $J$ is a complex structure on a Lie
algebra $\ggo$, then  the J-bracket $[\cdot , \cdot ]_J$  defined on $\ggo$ as follows:
\[ [x, y]_J = [Jx,y] + [x,Jy], \qquad  x, y \in \ggo,\]
satisfies the Jacobi identity. This follows from  the integrability condition for $J$ 
Moreover, $\ggo_J:=(\ggo, [\cdot, \cdot]_J)$ is a complex Lie
algebra. If $\ggo$ carries also an ad-invariant metric $g$ such
that $(J,g)$ is a Hermitian structure, then $\ggo$ admits a Lie
bialgebra structure. Such a Lie bialgebra is called of {\em
complex type} (see \cite{ABO} for more details). Since an even dimensional compact Lie algebra $\ggo$ admits a complex structure \cite{S} the modified cotangent of $\ggo$, $\nn(\ggo, \ad^{\ggo})$ can be equipped with a classical r-matrix. A more general case of r-matrices and Lie bialgebra structures on real simple Lie algebras was studied in \cite{AJ}.

\subsection{Constructing Lie algebras with ad-invariant metrics} 
Lie algebras provided with an ad-invariant metric can be obtained in
 the following way.
Let $(\bb, \varphi)$ be a  Lie algebra with ad-invariant metric $\varphi$ and let $S$ be a skew symmetric derivation of $(\bb, \varphi)$.  Consider the vector space direct sum $\RR Z \oplus \bb \oplus \RR T$ and equipp this vector space with the following Lie bracket: 
$$[z_1 Z + B_1+t_1 T, z_2 Z + B_2 + t_2 T]=\varphi(S B_1, B_2) Z + [B_1, B_2]_{\bb} + t_1 S B_2 - t_2 S B_1$$
where $z_i, t_i\in \RR, i=1,2,$ $B_1, B_2 \in \bb.$
 The metric $\la \,,\, \ra$ on $\ggo = \RR Z \oplus \bb \oplus \RR T$ obtained as a orthogonal extension of $\varphi$ by $\la \,, \,\ra_{\bb \times \bb}= \varphi$, $\la Z, T\ra=1$ and the other relations zero, allows to extend $\varphi$ to  an ad-invariant metric on $\ggo$. The Lie algebra $(\ggo, \la \,, \,\ra)$ is called the {\it double extension} of $(\bb, \varphi)$ by $(\RR, S)$.

It can be proved that any solvable Lie algebra $\ggo$  endowed with an ad-invariant  metric $\la \,, \,\ra$  is a double extension of a solvable Lie algebra with an ad-invariant metric $(\bb, \varphi)$ by  $\RR$ with a  certain skew symmetric derivation $S$ (see \cite{MR} for instance). 

The construction above shows that if $\ggo$ is 2-step nilpotent then $\bb$ should be a 2-step nilpotent Lie algebra and the fixed  skew symmetric derivation should be 2-step nilpotent. Our next result is related to the skew symmetric derivations of the Lie algebras $\nn(\vv, \rho)$. 

\begin{prop} Let $\nn(\vv, \rho)$ be a 2-step nilpotent Lie algebra of corank zero, with ad-invariant metric as in Theorem (\ref{t1}). Any skew symmetric derivation $\psi$ admits a matricial representation in a basis $\{z_i, v_j\}_{i,j=1}^n$ such such that $\la z_i, v_j \ra=\delta_{ij}$ as:
$$\left( \begin{matrix} -B^t & C \\
0 & B \end{matrix} \right) \qquad B\in M(n,\RR) : -\rho(B^t w) = \rho(w) B + B^t\rho(w)\quad  C \mbox{ skew-symmetric }
$$
\end{prop}
 The proof follows by canonical computations by asking a authomorphism to be skew symmetric in the desired basis.
 
 \vskip 4pt
 
 In the case of 2-step nilpotent Lie algebras any adjoint map $\ad_x$  will be represented by a matrix $C$ as above for any $v\in \vv$ (see (\ref{coef})). In particular if $\nn(\vv,\rho)$ is isomorphic to the the six dimensional 2-step nilpotent Lie algebra of corank zero, then from the extension  procedure we cannot get a 2-step nilpotent Lie algebra of corank zero, since we proved that there is no 2 step-nilpotent Lie algebra of corank zero admitting an ad-invariant metric with four dimensional commutator. Explicitely any skew symmetric derivation in this case will be representated by matrices $B$, $C$  as in the Proposition where  $B$ is symmetric of zero trace and $b_{11}=0$; the matrix  $C$ is  the matricial representation of an adjoint map.

\section{On 2-step nilpotent Lie groups provided with a bi-invariant metric} Our aim is to get information about the Lie groups that admit an ad-invariant metric. In the first part we work out these Lie groups making use of models provided by Riemannian geometry. In the second part we point out some features related to the geometry of a 2-step nilpotent Lie group when it is  endowed with a bi-invariant metric.

\subsection{A comparaison with the definite case}

People interested in the Riemannian geometry of a 2-step nilpotent Lie group equipped with a left invariant metric use to describe the  structure of the corresponding  Lie algebras in terms of certain maps $J_z$ with $z \in \zz$ constructed as follows. Let  $\nn$ be a 2-step nilpotent Lie algebra with an inner product $(\, , \,)$. Let $\vv$ denote the orthogonal complement of $\zz$ in $\nn$. Each element of $z \in \zz$ defines a skew symmetric linear map $J_z:\vv \to \vv$ given by $(J_z u, v) = (z, [u,v])$  for all $u,v\in \vv$. And one can see that the converse is true (see for instance \cite{E1}).

Next we seek  a characterization of 2-step Lie algebras which can be endowed with an ad-invariant metric making use  of the maps $J_z$ constructed up to a inner product.  In view of Theorem (\ref{t1}) we have

\vskip 4pt

{\bf Theorem  \ref{t1}' }{\it   Let $(\nn, ( \,, \, ))$ be a 2-step nilpotent Lie algebra with a inner product $(\, , \, )$. Then $\nn=\zz \oplus \vv$, where $\vv = \zz^{\perp}$,  can be endowed with an ad-invariant metric if and only if 

\begin{enumerate}
\item[\ri]  $C^1(\nn)$ is isomorphic to the orthogonal complement $\vv$;
\item[\rii]  $J:C^1(\nn) \to \End(\vv)$ is injective and 
\item[\riii] there exists a linear isomorphism $S: \vv \to C^1(\nn)$ satisfying $J_{Su} v + J_{Sv}u =0$ for all $u,v \in \vv$.
\end{enumerate}
} 

\begin{proof} If $\nn$ can be equipped with an ad-invariant metric then by Proposition (\ref{p1}) it admits a splitting $\nn= \tilde{\zz} \oplus \tilde{\nn}$ where $\tilde{\nn}$ is a 2m-dimensional 2-step nilpotent non degenerate ideal of corank zero containing the m-dimensional commutator $C^1(\nn)$. Thus $\dim \nn= 2m +k$ where $k=\dim \tilde{\zz}$. In fact the center $\zz$ of $\nn$ decomposes into a direct sum of vector subspaces $\zz = \tilde{\zz} \oplus C^1(\nn)$. If $\vv$ is a orthogonal complement of the center with respect to $( \,, \,)$ then its dimension must be  m and we have proved the first assertion. 

The map $J:C^1(\nn)  \to \End(\vv)$ is injective. In fact if $J_zu=0$, for all $u\in \vv$ then  $0=(J_zu, v)=(z, [u,v])$ for all $u,v \in \vv$. Hence since the inner product induces a isomorphim $C^1(\nn) \simeq \vv^{\ast} \simeq \vv$ we get  $z=0$ and we proved ii).

By  (\ref{coad}) there exists a linear isomorphism $T:\ggo \to \ggo^{\ast}$ satisfying $\ad^{\ast}_u = T \circ \ad_u \circ T^{-1}$ for all $u \in \vv$. Let $\ell:\ggo \to \ggo^{\ast}$ denote the isomorphism induced by the inner product, that is $\ell(x) (y) = ( x, y)$. We have 
$$\ad^{\ast}_u \ell(z) = -\ell(z) \circ \ad_u = -\ell(\ad^t_u z) = -\ell(J_z u)$$
where $\ad^t_u$ denotes the transpose of $\ad_u$ relative to the inner product $(\, , \,)$. Thus $-\ell(J_z u) = T \circ \ad_u \circ T^{-1} \ell(z)$ for all $z \in C^1(\nn)$, $u\in \vv$.  

Let $\tilde{T}:\nn \to \nn$ be the linear isomorphism given by $\tilde{T}= T^{-1} \circ \ell$ and denote  $\pi_{\vv}: \nn \to \vv$  the orthogonal projection over $\vv$. Setting $\tilde{S}= \pi_{\vv} \circ \tilde{T}$ it holds $\ad_u \tilde{T}z = \ad_u(\tilde{S}z)$ for all $z \in C^1(\nn)$.  Moreover $\tilde{S}:C^1(\nn) \to \vv$ is an isomorphism. In fact assume $\tilde{S} z=0$. Then we have $\ell(J_z u) = \tilde{T} ([u, \tilde{S}z]) = 0$ for all $u\in \vv$. Therefore $0=(J_z u, v)= (z, [u,v])$ for all $u,v\in \vv$, implying $z=0$. Let $S= \tilde{S}^{-1}:\vv \to C^1(\nn)$. Since $J_{Su} u = -\tilde{T} \ad_u u=0$ we proved the third assertion.
 
For the converse notice that $\rho= J_S:\vv \to \sog(\vv)$ satisfies requirements of Theorem (\ref{t1}) in the 2-step nilpotent ideal $\tilde{\nn}= C^1(\nn) \oplus \vv$, and this  allows to say that the canonical metric $\la \,, \,\ra$ on $(\vv^{\ast} \simeq C^1(\nn))\oplus \vv$ defines  an ad-invariant metric on $\tilde{\nn}$ as in (\ref{hyp})  with the Lie bracket given in (\ref{t1}), where  $\la u, [v,w]\ra= (J_{Su} v, w)$.
By extending the hyperbolic metric of $C^1(\nn) \oplus \vv$ to a metric on $\zz/ C^1(\nn)$ we get an ad-invariant  metric on $\nn$.
\end{proof}

 %The previous theorem shows a necessary condition for a 2-step nilpotent Lie algebra to admit an ad-invariant metric, as we point out in the following corollary. 
 
 A 2-step nilpotent Lie algebra $\nn$ is non singular if $\ad_x: \nn \to \zz$  is surjective for all $x \in \nn - \zz$. Examples  are the 2-step nilpotent Lie algebras of H-type. 
 Recall that a  2-step nilpotent Lie algebra equipped with a inner product $(\, , \,)$ is said to be of {\em H-type} if for  all $z\in \zz$ the maps $J_z$ satisfy $J_z^2= -(z,z)Id,$ where $Id$ denotes the identity on $\vv=\zz^{\perp}$. 

\begin{cor} A non singular 2-step nilpotent Lie algebra does not admit an ad-invariant metric. In particular 2-step nilpotent Lie algebras of H-type cannot be equipped with an ad-invariant metric.
\end{cor} 
\begin{proof}It is known that a  2-step nilpotent Lie algebra is non singular if and only if  for any inner product $(\, , \,)$ on $\nn$ the corresponding maps $J_z$ where $z \in \zz$ are non singular on $\vv = \zz^{\perp}$ (see for instance \cite{E1}). Let $\nn$ be a 2-step nilpotent Lie algebra admitting an ad-invariant metric. By
the Theorem above (\ref{t1}')  for every non zero $z\in C^1(\nn)$ there is a non zero $u \in \vv$ such that $J_z u=0$, thus the kernel of $J_z$ is non trivial.  This implies that $\nn$ cannot be non singular. 

The second assertion follows from the first one since any $H$-type Lie algebra is nonsingular.
\end{proof}

\begin{lem} Let $(\nn, (\,, \,))$ be a 2-step nilpotent Lie algebra of corank zero with inner product $(\, , \,)$  and $\dim z = \frac 12 \nn$. Assume $J:\zz \to \sog(\vv)$ is injective and that  there is a basis $\{z_1, z_2, \hdots, z_n\}$ of $\zz$ for which there exists  a basis $\{v_1, v_2, \hdots, v_n\}$ of $\vv$ satisfying: $v_i \in Ker J_{z_i}$  for all i and 
\begin{equation}\label{ev}
\mbox{ for } v_i \in Ker J_{z_i}, v_j \in Ker J_{z_j} \quad \mbox{ it holds }\quad  v_i + v_j \in  Ker J_{z_i + z_j}.
\end{equation}
Then $\nn$ can be endowed with an ad-invariant metric.
\end{lem}
\begin{proof} One proves that the condition above allows to define $S: \vv \to C^1(\nn)$ by $S(v_i) = z_i$, which satisfies  $J_{Su} v + J_{Sv}u =0$ for all $u,v \in \vv$.
\end{proof}

\subsection{On the geometry arising from a bi-invariant metric}
Let $G$ be  a connected Lie group with Lie algebra $\ggo$ provided with a metric $\la \,, \,\ra$ 
  and let $Q$ be the left invariant pseudo-Riemannian metric on $G$ obtained by left translations such that $Q_e = \la \,, \,\ra$. The fact that $Q$ is bi-invariant (i.e. also right invariant) is equivalent to either of the following statements:
 
 i) $\la \,, \,\ra $ satisfies (\ref{ad});
 
 ii) the geodesics through the identity $e$ are the one-parameter subgroups of $G$.

 Moreover if $\nabla$ denotes the Levi Civita connection for $Q$, $R$ the curvature tensor, $Ric$  the Ricci tensor and $B$ the Killing form on $\ggo$ one has the formulas
 
 $$\begin{array}{rcl}
 \nabla_x y & = & \frac12 [x, y] \\
 R(x,y) & = & - \frac14 \ad([x,y]) \\
 Ric & = & -\frac14 B
 \end{array}
 $$
 for all left invariant vector fields $x,y \in \ggo$.  For the Lie algebra $\mathfrak{hol}(G,Q)$ of the holonomy group $(G,Q)$ we have
 $$\mathfrak{hol}(G,Q) = span\{ R(x,y)\, /\,  x, y \in \ggo\} = \ad[\ggo, \ggo] \subset \sog(\ggo, \la , \ra)
 $$
  Thus a direct application of these formulas to the case of 2-step nilpotent Lie groups $N$ endowed with a bi-invariant metric $Q$ proves the following proposition.

 \begin{prop} Any  simply connected 2-step nilpotent Lie group $(N, Q)$ attached with a bi-invariant metric $Q$ is flat (therefore with trivial  holonomy). 
 \end{prop}
 
 Our next goal is to study the group of isometries of $(N,Q)$.   It is known that this group denoted $I(N)$  endowed with the compact-open topology, is a Lie group, even a Lie transformation group of $N$. Let $F(N)$ denote the closed subgroup of $I(N)$ consisting of those isometries of $N$ fixing the identity element of $N$ and let $L(N)$ denote the subgroup of left translations $L(N)$. Then $I(N)=L(N)F(N)$. 
 
 In \cite{Mu} a theory  about the isometries of a connected Lie group provided with a bi-invariant metric was established. It is proved that $L(N)$ is a closed connected subgroup isomorphic to $N$, such that $L(N)\cap F(N)=\{id\}$ being $id$ the identity map  of $N$. Since every local isometry at $e$ can be extended (uniquely) to a global isometry if $N$ is simply connected, then any element of  $F(N)$ is uniquely determined by its differential at $e$, $d\Phi_e$.
 
 \begin{thm}\cite{Mu} \label{mu} Let $G$ be  a connected Lie group with Lie algebra $\ggo$ provided with a an ad-invariant metric $\la \,, \,\ra$ 
  and let $Q$ be the left invariant pseudo-Riemannian metric on $G$ obtained by left translations such that $Q_e = \la \,, \,\ra$.  Let $A$ be a linear endomorphism of  $\ggo$. Then there exists a local isometry $\Phi$ of $G$ at $e$ with $d\Phi_e=A$ if and only if $A$ satisfies the following conditions:
 
 \begin{enumerate}
 \item[\ri] $\la Ax, Ay\ra = \la x,y\ra$ for all $x,y \in \ggo$.
 
 \item[\rii] $A([x, [x,y]]) = [Ax, [Ax, Ay]]$ for all $x,y \in \ggo$.
 \end{enumerate}
 \end{thm}
 
 We shall say that the corank of the 2-step nilpotent Lie group $N$ is $k$ if this is the corank of its Lie algebra $\nn$.
 If $N$ is a simply connected 2-step nilpotent Lie group of corank zero with bi-invariant metric $Q$ the previous result says that the group of isometries is like in the abelian case. In fact, a  direct application of the Theorem above implies that the  group of isometries fixing the identity  is $O(k,k)$ being  $k$ the dimension of the commutator.  
 Hence the group of isometries of $N$ can be much larger than the group of all automorphism which leave the metric invariant (see (\ref{t2})).

 \begin{thm} Let $(N, Q)$ be a  connected 2-step nilpotent Lie group of corank zero provided with a bi-invariant metric $Q$. Then the isometry group is $I(N)= O(k,k)N$, where $k$ is the dimension of the commutator subgroup of $N$.
 \end{thm}

Let $N$, $N'$  denote the simply connected 2-step nilpotent Lie groups of corank zero with respective Lie algebras $\nn(\vv,\rho)$, $\nn(\vv,\rho')$ as in  Theorem (\ref{t1}) constructed from a vector space $(\vv , \la \, , \,\ra_{\vv})$ with $\rho$ and $\rho'$ respectively (Theorem (\ref{t1})).  The Theorem of M\"uller proves the following result.

\begin{thm}  Let $N$ and $N'$ be simply connected  2-step nilpotent Lie groups of corank zero of the same dimension  endowed with a bi-invariant pseudo-Riemannian metric as in the previous paragraph. Then $N$ and $N'$ are isometric.
\end{thm}
\begin{proof} In view of Theorem (\ref{mu}) any isometry fixing the identity element does not depend on the algebraic structure therefore $N$ and $N'$ are isometric if and only if their Lie algebras are isomorphic as vector spaces  and there exists a map $A: \nn \to \nn'$ satisfying (i). Since $\nn$ and $\nn'$ have the same dimension assume this is the even integer $2m$. Let  $\{z_i, v_j\}_{i,j=1}^m$, $\{z_i', v_j'\}_{i,j=1}^m$  be  basis of $\nn$ and $\nn'$ respectively satisfying $\la z_i, v_j\ra = \delta_{ij} = \la v_i', z_j'\ra $ and $0=\la z_i, z_j\ra =  \la v_i, v_j\ra = \la z_i', z_j'\ra = \la v_i', v_j'\ra$ for all $i,j= 1, \hdots, m$. The existence of such basis follows from Theorem (\ref{t1}) and Lemma (\ref{l1}). 
Let $A: \nn \to \nn'$ be the linear map defined by $A z_i=z_i'$ and $A v_i=v_i'$ for $i=1, \hdots, m$. Then $A$ is a orthogonal isomorphism of the vector spaces $\nn$ and $\nn'$, satisfying (i) and it  extends to an isometry from $N$ to $N'$.
\end{proof}

%\begin{rem} Notice that the assumption of zero corank in the previous Theorem cannot be removed. In fact  it is possible to have isometric 2-step nilpotent Lie groups of the same dimension but with different coranks. 
%\end{rem}

 \end{document}